\documentclass[reqno,12pt]{amsart}
\usepackage{amsmath, latexsym, amsfonts, amssymb, amsthm, amscd}

\numberwithin{equation}{section}

\newtheorem{theorem}{Theorem}[section]
\newtheorem{lemma}[theorem]{Lemma}
\newtheorem{prop}[theorem]{Proposition}
\newtheorem{cor}[theorem]{Corollary}
\newtheorem{rem}[theorem]{Remark}

\newcommand{\gga}{\gamma}            

\newcommand{\gep}{\varepsilon}       

\newcommand{\cB}{{\ensuremath{\mathcal B}} }
\newcommand{\cF}{{\ensuremath{\mathcal F}} }
\newcommand{\cP}{{\ensuremath{\mathcal P}} }
\newcommand{\cE}{{\ensuremath{\mathcal E}} }

\newcommand{\cL}{{\ensuremath{\mathcal L}} }

\newcommand{\cD}{{\ensuremath{\mathcal D}} }

\newcommand{\cI}{{\ensuremath{\mathcal I}} }

\newcommand{\E}{{\ensuremath{\mathbb E}} }

\newcommand{\bbP}{{\ensuremath{\mathbb P}} }

\newcommand{\R}{{\ensuremath{\mathbb R}} }

\date{}

\title{Quasi-invariance properties of a class of subordinators}

\author{Max-K. von Renesse}
\address{ Institut f\"ur Mathematik,
TU Berlin,  Strasse des 17. Juni 136, 10623 Berlin, Germany }
\email{mrenesse\@@math.tu-berlin.de}

\author{Marc Yor}
\address{Laboratoire de Probabilit{\'e}s et Mod\`eles Al\'eatoires (CNRS U.M.R. 7599) \\  Universit{\'e} Paris 6
-- Pierre et Marie Curie, U.F.R. Mathematiques, Case 188, 4 place
Jussieu, 75252 Paris cedex 05, France }

\author{Lorenzo Zambotti}
\address{Laboratoire de Probabilit{\'e}s et Mod\`eles Al\'eatoires (CNRS U.M.R. 7599) \\ Universit{\'e} Paris 6
-- Pierre et Marie Curie, U.F.R. Mathematiques, Case 188, 4 place
Jussieu, 75252 Paris cedex 05, France }
\email{zambotti\@@ccr.jussieu.fr}

\begin{document}

\begin{abstract}
We study absolute-continuity properties of a class of stochastic
processes, including the gamma and the Dirichlet processes. We prove
that the laws of a general class of non-linear transformations of
such processes are locally equivalent to the law of the original
process and we compute explicitly the associated Radon-Nikodym
densities. This work unifies and generalizes to random non-linear
transformations several previous results on quasi-invariance of
gamma and Dirichlet processes.
\end{abstract}


\maketitle

\section{Introduction}

In this paper we present several absolute-continuity results
concerning, among others, the gamma process and the Dirichlet
processes. We recall that the gamma process $(\gga_t)_{t\geq 0}$ is
a subordinator, i.e. a non-decreasing L\'evy process, with gamma
marginals, i.e. $\gamma_0=0$ and
\[
\bbP(\gga_t\in dx) = p_t(x) dx, \qquad p_t(x):= 1_{[0,\infty)}(x)\,
\frac1{\Gamma(t)} \, x^{t-1} \, e^{-x}, \qquad t>0, \ x\in\R.
\]
Moreover for any $T>0$, we define the Dirichlet process over $[0,T]$
as $D_t^{(T)}:=\gamma_t/\gamma_T$, $t\in[0,T]$; we recall that
$\gamma_T$ is independent of $(\gamma_t/\gamma_T, t\in[0,T])$ and
that, therefore, $(D_t^{(T)}, t\in[0,T])$ is equal in law to the
gamma process conditioned on $\{\gga_T=1\}$. See \cite{yor} for a
survey of the main properties of the gamma process.

The gamma process has been the object of intense research activity
in recent years, both from pure and applied perspectives, such as in
representation theory of infinite dimensional groups, in
mathematical finance and in mathematical biology (see e.g.
\cite{tvy, baxter,handa}). Quasi-invariance properties of the
associated probability measure on path or measure space with respect
to canonical transformations often play a central role. We recall
that, given a measure $\mu$ on a space $X$ and a measurable map
$T:X\mapsto X$, quasi-invariance of $\mu$ under $T$ means that $\mu$
and the image measure $T_*\mu$ are equivalent, i.e. mutually
absolutely continuous. A classical example is the Girsanov formula
for additive perturbations of Brownian motion (see, e.g.,
\cite{reyo}, Chap. VIII).

In this paper we study quasi-invariance properties for a class of
subordinators which we denote by $(\cL)$ and define below, with
respect to a large class of non-linear sample path transformations.
In particular, we unify and extend previous results on the real
valued gamma and Dirichlet processes.

Quasi-invariance properties of L\'evy processes have been studied
for quite some time, see e.g. Sato \cite[p. 217-218]{sato}. In the
case of the gamma process, for any measurable function
$a:\R_+\mapsto\R_+$ with $a$ and $1/a$ bounded, the laws of
$(\int_0^t a_s \, d\gamma_s, t\geq 0)$ and $(\gamma_t, t\geq 0)$ are
locally equivalent, see \cite{tvy}. By local equivalence of  two
real-valued processes $(\eta_t, t\geq 0)$ and $(\zeta_t, t\geq 0)$,
we mean that for all $T>0$ the laws of $(\eta_t, t\in[0,T])$ and
$(\zeta_t, t\in[0,T])$ are equivalent.

Here, we show the same property for a much wider class of
transformations $\xi\mapsto K(\xi)$, e.g. $(\xi_t)_{t\geq
0}\mapsto(K(t,\xi_t))_{t\geq 0}$ and $(\xi_t)_{t\geq
0}\mapsto(\sum_{s \leq t }K(s,\Delta\xi_s))_{t\geq 0}$, where
$(\xi_t)_{t\geq 0}$ is a $(\cL)$-subordinator
 and $K(s,.)$ is a $C^{1,\alpha}$-isomorphism of $\R_+$
for each $s\geq 0$ and $\alpha\in]0,1[$. Using the mentioned
properties of the gamma process, we establish analogous
quasi-invariance results for transformations $D^{(T)}\mapsto
K(D^{(T)})$ of the Dirichlet process, e.g.
$(D_t^{(T)})_{t\in[0,T]}\mapsto(K(t,D_t^{(T)}))_{t\in[0,T]}$,
where $K(s,.)$ is an increasing
$C^{1,\alpha}$-isomorphism of $[0,1]$ for each $s\in [0,T]$.

In all these cases, we compute the Radon-Nikodym density explicitly
and study its martingale structure. We notice that our approach
allows to treat the previously mentioned results  by
Vershik-Tsilevich-Yor \cite{tvy,tvy2}, together with Handa's
\cite{handa} and the recent one  by  Renesse-Sturm \cite{vrs} on
Dirichlet processes, within a unified framework.

The paper ends with an application to SDEs driven by
$(\cL)$-subordinators. Finally we point out that, in the same spirit
as in \cite{vrs}, each quasi-invariance property we show yields
easily an integration by parts formula on the path space; such
formulae can be used in order to study an appropriate Dirichlet form
and the associated infinite-dimensional diffusion process. These
applications will be developed in a future work.

\subsection{The main result}\label{setting}

Let $(\xi_t)_{t\geq 0}$ be a subordinator, i.e. an increasing L\'evy
process with $\xi_0=0$. In this paper we consider subordinators in
the class $(\cL)$, meaning {\it with logarithmic singularity}, i.e.
we assume that $\xi$ has zero drift and L\'evy measure
\[
\nu(dx) \, = \, g(x)\, dx, \qquad  x> 0,
\]
where $g:]0,\infty[\mapsto\R_+$ is measurable and satisfies
\begin{enumerate}
\item[(H1)] $g>0$ \ and \
$\displaystyle{\int_1^\infty g(x) \, dx < \infty}$;
\item[(H2)] there exist $g_0\geq 0$ and $\zeta:[0,1]\mapsto\R$ measurable such that
\[
g(x) = \frac{g_0}x + \zeta(x),  \qquad \forall x\in\, ]0,1],
\qquad \text{and}  \quad \int_0^1|\zeta(x)| \, dx < +\infty.
\]
\end{enumerate}
We recall that for all $t\geq 0$, $\lambda>0$
\[
\E\left( e^{-\lambda \xi_t} \right) =
\exp\left(-t\Psi(\lambda)\right), \qquad \Psi(\lambda) :=
\int_0^\infty \left(1-e^{-\lambda x}\right) \, g(x) \, dx.
\]
For the general theory of subordinators, see \cite{bertoin}. We
denote by $\cF_t:=\sigma(\xi_s: s\leq t)$, $t\geq 0$, the filtration
of $\xi$. We denote the space of c\`adl\`ag functions on $[0,t]$ by
$\cD([0,t])$, endowed with the Skorohod topology.

\begin{rem}{\rm
In the particular case of the gamma process $(\gga_t)_{t\geq 0}$,
mentioned above, we have
\[
g(x) = \frac{e^{-x}}x, \quad x>0, \qquad
\Psi(\lambda)=\log(1+\lambda), \quad \lambda\geq 0.
\]
}
\end{rem}

\noindent
We consider a measurable function
$h:\R_+\times\Omega\times\R_+\mapsto\R_+$ such that
\begin{enumerate}
\item $h$ is $\cP\otimes\cB_{\R_+}$-measurable, where $\cP$ denotes the
predictable $\sigma$-algebra generated by $\xi$;
\item denoting $h(s,a)=h(s,\omega,a)$, there exist finite constants
$\kappa>1$ and $\alpha\in]0,1[$, such that almost surely
\begin{equation}\label{H'1}
|h(s,x)-h(s,y)| \leq \kappa |x-y|^\alpha, \qquad \forall \,
x,y\in\R_+, \ s\geq 0,
\end{equation}
\begin{equation}\label{H'2}
0< \kappa^{-1} \leq h(s,x) \leq \kappa< \infty, \qquad \forall \,
x\in\R_+, \ s\geq 0.
\end{equation}
\end{enumerate}
Then we set
\[
H(s,x) = \int_0^x h(s,y) \, dy, \qquad \forall \, x\geq 0, \ s\geq
0.
\]
Notice that a.s. $H(s,\cdot):\R_+\mapsto\R_+$ is necessarily a
$C^1$-diffeomorphism for all $s\geq 0$. We set
\[
\Delta \xi_s := \xi_s-\xi_{s-}, \qquad s\geq 0,
\]
and for convenience of notation
\begin{equation}\label{notation}
h(s,0)\cdot \frac{g(H(s,0))}{g(0)} := 1, \qquad \forall \, s\geq 0.
\end{equation}
We can now state the main result of this paper
\begin{theorem}\label{main}$ $
\begin{enumerate}
\item The process
\[
M_t^H :=  \exp\left( g_0 \int_0^t \log h(s,0) \, ds \right)
\prod_{s\leq t} \left[h(s,\Delta \xi_s) \cdot \frac{g(H(s,\Delta
\xi_s))}{g(\Delta \xi_s)}\right], \quad t\geq 0,
\]
is a $(\cF_t,\bbP)$-martingale with $\E(M_t^H)=1$ and a.s.
$M_t^H>0$. We can uniquely define a probability measure $\bbP^H$
such that $\bbP^H_{|\cF_t}=M_t^H\cdot \bbP_{|\cF_t}$ for all $t\geq
0$.
\item Setting
\begin{equation}\label{defxih}
\xi^H_t := \sum_{s\leq t} H(s,\Delta \xi_s), \qquad t\geq 0,
\end{equation}
then $\xi^H$ is distributed under $\bbP^H$ as $\xi$ under $\bbP$.
\end{enumerate}
\end{theorem}
Note that Theorem \ref{main} is a local equivalence result for the
laws of $\xi$ and $\xi^H$, since a.s. $M_t^H>0$. The theorem is
stated for general subordinators in the class $(\cL)$ defined above
and for a general random transformation; in section
\ref{section:special} we consider some special cases of the general
result, and in section \ref{section:dirichlet} we consider the case
of the Dirichlet process.

\subsection{A parallel between the gamma process and Brownian
motion} \label{parallel} The absolute-continuity results presented
in this paper can be better understood by comparison with some
analogous properties of Brownian motion.

The Girsanov theorem for a Brownian motion $(B_t, t\geq 0)$ states the
following property: if $(a_s, s\geq 0)$ is an adapted and (say)
bounded process, then the law of the process
\[
t\mapsto B_t + \int_0^t a_s \, ds, \qquad t\geq 0,
\]
is locally equivalent to that of $(B_t, t\geq 0)$, with explicit
Radon-Nikodym density. We call this property quasi-invariance by
addition.

As a byproduct case of our Theorem \ref{main}, the gamma process
$\gamma$ has an analogous property of quasi-invariance by
multiplication (see also \cite{tvy}): if $(a_s, s\geq 0)$ is a
predictable process such that $a$ and $1/a$ are bounded, then the
law of
\[
t\mapsto \int_0^t a_s \, d\gamma_s, \qquad t\geq 0,
\]
is locally equivalent to that of $(\gamma_t, t\geq 0)$, and we
compute explicitly the Radon-Nikodym density. In fact, we can prove
the same quasi-invariance property for all $(\cL)$-subordinators.

The Girsanov theorem for Brownian motion has important applications
in the study of stochastic differential equations (SDEs) driven by a
Wiener process; likewise, our Theorem \ref{main} allows to give
analogous applications to SDEs driven by $(\cL)$-subordinators, e.g.
to compute explicitly laws of solutions; see section \ref{sdes}.

\section{A generalization of a formula of Tsilevich-Vershik-Yor}

Within the framework of subsection \ref{setting}, the law of $\xi$
with $\xi_0=0$ is characterized by its Laplace transform, i.e. for
any measurable bounded $\lambda:\R_+\mapsto\R_+$
\[
\E\left[\exp\left(-\int_0^t \lambda_s \, d\xi_s\right)\right] =
\exp\left(-\int_0^t \Psi(\lambda_s)\, ds\right).
\]
In order to prove Theorem \ref{main}, we shall show that $\xi^H$
has, under $\bbP^H$, the same Laplace transform as $\xi$ under
$\bbP$, namely for all measurable bounded $\lambda:\R_+\mapsto\R_+$
\[
\E^H\left[\exp\left(-\int_0^t \lambda_s \, d\xi^H_s\right)\right] =
\exp\left(-\int_0^t \Psi(\lambda_s)\, ds\right).
\]
To do that, we shall show that the process
\[
\exp\left(-\int_0^t \lambda_s \, d\xi^H_s + \int_0^t
\Psi(\lambda_s)\, ds\right), \qquad t\geq 0,
\]
is a $((\cF_t),\bbP^H)$-martingale, which is equivalent to prove the
following
\begin{prop}\label{prmain}
 We set for all $t\geq 0$
\begin{equation}\label{defMH}
\begin{split}
M_t^{H,\lambda} :=  & \exp\left( \int_0^t \left( g_0  \log h(s,0) +
\Psi(\lambda_s)\right)\, ds \right) \cdot
\\ & \cdot \prod_{s\leq t} \left[
h(s,\Delta \xi_s) \, \frac{g(H(s,\Delta \xi_s))}{g(\Delta \xi_s)} \,
\exp\left( -\lambda_{s}  H(s,\Delta \xi_s)\right)\right].
\end{split}
\end{equation}
Then $M^{H,\lambda}$ is a $(\cF_t,\bbP)$-martingale with
$\E(M_t^{H,\lambda})=1$ and a.s. $M_t^{H,\lambda}>0$.
\end{prop}
\noindent Tsilevich-Vershik-Yor prove in \cite{tvy} the same result
for $\xi$ a gamma process and $H(s,x)=c(s)x$, for
$c:\R_+\mapsto\R_+$ measurable and deterministic.

We say that a real-valued process $(\zeta_t, t\geq 0)$ has bounded variation,
if a.s. for all $T>0$ the real-valued function $[0,T]\ni t\mapsto \zeta_t$ has
bounded variation.

\medskip
\begin{lemma}\label{pr1}
Let $F:\R_+\times\Omega\times\R_+\mapsto ]-1,\infty[$ such that
\begin{itemize}
\item $F$ is $\cP\otimes\cB_{\R_+}$-measurable, where $\cP$ denotes the
predictable $\sigma$-algebra generated by $\xi$;
\item there exists a finite constant $C_F$ such that a.s. for almost every $s\geq
0$
\begin{equation}\label{stima}
F(s,0)=0, \qquad \int \nu(dx) \, \E \left[\Big| F(s,x) \Big| \right]
\, \leq \, C_F < \infty.
\end{equation}
\end{itemize}
Then
\begin{enumerate}
\item the process
\[
x^F_t := \sum_{s\leq t} F(s,\Delta\xi_s) - \int_0^t ds \int \nu(dx)
\, F(s,x), \qquad t\geq 0
\]
is a martingale with bounded variation;
\item the process
\[
\cE^F_t := \exp\left(-\int_0^t ds \int \nu(dx) \, F(s,x) \right) \,
\prod_{s\leq t} \Big( 1+F(s, \Delta\xi_s) \Big)
\]
satisfies
\begin{equation}\label{doleans}
\cE^F_t = 1 + \int_0^t \cE^F_{s-} \, dx^F_s, \qquad t\geq 0.
\end{equation}
Moreover $(\cE^F_t)$ is a martingale with bounded variation which
satisfies
\[
\E\left( \int_0^t \left| d\cE^F_u \right| \right) \leq 2 \, C_F \,
t.
\]
\item for all $t\geq 0$, a.s. $\cE^F_t>0$.
\end{enumerate}
\end{lemma}
\noindent{\it Proof}. Notice first that $x^F$ is well defined, since
by \eqref{stima}
\[
\E\left[ \sum_{s\leq t} \left| F(s,\Delta\xi_s) \right| \right] =
\int_0^t ds \int \nu(dx) \, \E \left[\Big| F(s,x) \Big| \right] \,
\leq \, C_F \, t < \infty.
\]
Since $\{(s,\Delta\xi_s), s\geq 0\}$ is a Poisson point process
with intensity measure $ds \, \nu(dx)$, it follows
immediately that $x^F$ is a local martingale. Furthermore, a.s. the
paths of $x^F$ have bounded variation, since
\[
\E\left( \int_0^t \left| dx^F_s \right| \right) \leq 2 \, C_F \, t,
\qquad t\geq 0.
\]
Therefore, $(x^F_t,t\geq 0)$ is a true martingale; indeed, for any
$t>0$, $\sup_{s\leq t} |x^F_s|\leq \int_0^t \left| dx^F_s \right|$,
and therefore
\[
\E\left( \sup_{s\leq t} |x^F_s| \right) \leq 2 \, C_F \, t, \qquad
t\geq 0;
\]
by Proposition IV.1.7 of \cite{reyo} we obtain the claim.

Since $\cE^F$ is the Dol\'eans exponential associated with the
martingale $x^F$, i.e. it satisfies \eqref{doleans}, it is clear
that $\cE^F$ is a local martingale (see chapter 5 of \cite {apple}).
Moreover, since $\cE^F$ is non-negative, then it is a
super-martingale and in particular $\E\left( \cE^F_t \right)\leq
\E\left( \cE^F_0 \right)=1$. Furthermore,
\[
\E\left( \int_0^t \left| d\cE^F_u \right| \right) = \E\left(
\int_0^t \cE^F_{u-} \, \left| dx^F_u \right| \right) \leq \int_0^t
\E\left( \cE^F_u \right) 2 \, C_F \, du \leq 2 \, C_F \, t.
\]
The same argument as for $x^F$ yields:
\[
\E\left( \sup_{s\leq t} \left| \cE^F_s \right| \right)  \leq 2 \,
C_F \, t, \qquad t\geq0,
\]
and therefore $\cE^F$ is a martingale.

In order to prove that $\cE^F_t>0$ a.s., by \eqref{stima} it is
enough to show that
\[
\log \prod_{s\leq t} \Big( 1+F(s, \Delta\xi_s) \Big) = \sum_{s\leq
t}  \log\Big( 1+F(s, \Delta\xi_s) \Big)> -\infty.
\]
Since $F(s, \Delta\xi_s)=\Delta x^F_s=x^F_s-x^F_{s-}>-1$, and $x^F$
has a.s. bounded variation, then there is a.s. only a finite number
of $s\in[0,t]$ such that $\Delta x^F_s<-1/2$ and therefore a.s.
$\inf_{s\leq t}\Delta x^F_s=:C_t>-1$. It follows that
\[
\sum_{s\leq t}  \log\Big( 1+\Delta x^F_s \Big) \geq
-\frac 1{C_t+1} \sum_{s\leq t}  \left| \Delta x^F_s\right| = -\frac
1{C_t+1}\int_0^t \left| dx^F_u \right| > - \infty, \quad \text{a.s.}
\qed
\]

\medskip\noindent
The main steps in the proofs of Proposition \ref{prmain} and Theorem
\ref{main} are the estimate \eqref{c_F} and the identity
\eqref{identity} below, which allow to apply Lemma \ref{pr1} to
\begin{equation}\label{defF} F(s,0):=0, \qquad
F(s,x) := h(s,x) \cdot \frac{g(H(s,x))}{g(x)} \cdot e^{-\lambda_s
H(s,x)}- 1, \quad x>0.
\end{equation}
\begin{lemma}\label{est0}
Let $\phi:\R_+ \mapsto\R_+$ a $C^1$ function such that $\phi(0)=0$,
\[
0< \kappa^{-1} \leq \phi'(x) \leq \kappa< \infty,
\qquad |\phi'(x)-\phi'(y)| \leq \kappa |x-y|^\alpha, \qquad \forall \,
x,y\in\R_+,
\]
where $\kappa>0$ and $\alpha\in]0,1[$. We set for all $a\geq 0$
\[
F_{a,\phi}: (0,\infty)\mapsto\R, \qquad F_{a,\phi} := \phi' \cdot
\frac{g(\phi)}{g} \cdot e^{-a\phi}- 1.
\]
There exists a finite constant $C=C(\kappa,\alpha,a)$ such that
\begin{equation}\label{c_F}
\int_0^\infty \Big| F_{a,\phi}(x) \Big|\, g(x)\, dx \leq C,
\end{equation}
and
\begin{equation}\label{identity}
\int_0^\infty  F_{a,\phi}(x) \, g(x)\, dx = -\Psi(a) -g_0 \log
\phi'(0).
\end{equation}
\end{lemma}
\noindent{\it Proof}. Notice that $\phi:\R_+\mapsto\R_+$ is a
diffeomorphism. First we have
\[
\begin{split}
& \int_{\kappa^{-1}}^\infty \Big| F_{a,\phi} \Big|\, g \, dx \leq
\int_{\kappa^{-1}}^\infty \phi' \, g(\phi) \, dx +
\int_{\kappa^{-1}}^\infty g\, dx \\ & =
\int_{\phi(\kappa^{-1})}^\infty g(y)\,  dy +
\int_{\kappa^{-1}}^\infty g\, dx \leq 2\int_{\kappa^{-2}}^\infty
g(x)\, dx < \infty.
\end{split}
\]
Now
\[
\begin{split}
& \int_0^{\kappa^{-1}} \Big| F_{a,\phi} \Big|\, g \, dx =
\int_0^{\kappa^{-1}} \Big| \phi' \, g(\phi) \, e^{-a \phi} -
g\Big|\, dx  \\ & \leq \int_0^{\kappa^{-1}} \phi' \, g(\phi) \left(
1- e^{-a \phi}\right) dx + \int_0^{\kappa^{-1}} \phi' \left|
g(\phi)-\frac{g_0}{\phi} \right| \, dx + \int_0^1 g_0 \left|
\frac{\phi'(x)}{\phi(x)} - \frac 1x\right| \, dx \\ & +
\int_0^{\kappa^{-1}} \left| \frac{g_0}x- g(x)\right|\, dx +
\int_0^{\kappa^{-1}} g(x) \, (1-e^{-a x}) \, dx =: I_0 + I_1 + I_2 +
I_3 + I_4.
\end{split}
\]
First we estimate $I_2$.
\[
\begin{split}
I_2 & =  \int_0^1 g_0 \left| \frac{\phi'(x)}{\phi(x)} - \frac
1x\right| \, dx = g_0\int_0^1 \left|
\frac{\phi(x)-x\phi'(x)}{x\phi(x)} \right| \, dx \\ & \leq \int_0^1
\frac{g_0}{\kappa^{-1}x^2} \left| \int_0^x
\left[\phi'(y)-\phi'(x)\right] dy \right| \, dx \leq g_0 \, \kappa^2
\int_0^1 \frac1{x^2} \int_0^x y^\alpha dy \, dx \leq \frac{g_0 \,
\kappa^2}{\alpha(1+\alpha)}.
\end{split}
\]
Recall now that $g(x)=\frac{g_0}x+\zeta(x)$ by (H2) above. Then
$I_3$ and $I_4$ can be estimated by
\[
I_3 = \int_0^{\kappa^{-1}}  \left| \frac{g_0}x- g(x)\right|\, dx
\leq \int_0^1 |\zeta|\, dx,
\]
and
\[
I_4  = \int_0^{\kappa^{-1}} g(x) \, (1-e^{-a x}) \, dx \leq \int_0^1
g_0 \, ax \, dx + \int_0^1 |\zeta|\, dx \leq ag_0 + \int_0^1
|\zeta|\, dx.
\]
Then $I_0$ and $I_1$ can be estimated similarly by changing variable
\[
I_1 = \int_0^{\kappa^{-1}} \phi' \left| g(\phi)-\frac{g_0}{\phi}
\right| \, dx = \int_0^{\phi(\kappa^{-1})} \left| g(x)- \frac{g_0}x
\right|\, dx \leq \int_0^1 |\zeta|\, dx,
\]
and
\[
I_0  = \int_0^{\kappa^{-1}} \phi' \, g(\phi) \left( 1- e^{-a
\phi}\right) dx = \int_0^{\phi(\kappa^{-1})} g(x) \, (1-e^{-a x}) \,
dx \leq ag_0 + \int_0^1 |\zeta|\, dx,
\]
since $\phi(\kappa^{-1})\leq 1$. Therefore, we have obtained
\[
\int_0^\infty \Big| F_{a,\phi} \Big|\, g \, dx \leq 2ag_0 +
\frac{g_0 \, \kappa^2}{\alpha(1+\alpha)} +
2\int_{\kappa^{-1}}^\infty g(y) \, dy + 4\int_0^1 |\zeta(x)|\, dx,
\]
and \eqref{c_F} is proven.

We turn now to the proof of \eqref{identity}. By \eqref{c_F} and
dominated convergence
\[
\int_0^\infty  F_{a,\phi} \, g\, dx = \lim_{\gep\searrow 0}
\int_\gep^\infty  F_{a,\phi} \, g\, dx.
\]
For all $\gep>0$ we have
\[
\int_\gep^\infty \phi' \, g(\phi) \, e^{-a \phi} \, dx = \Big[
y=\phi(x) \Big] = \int_{\phi(\gep)}^\infty g(y) \, e^{-a y} \, dy.
\]
Then we want to compute the limit as $\gep\searrow 0$ of
\[
\begin{split}
& \int_\gep^\infty  F_{a,\phi} \, g\, dx = \int_{\phi(\gep)}^\infty
g(x) \, e^{-a x} \, dx - \int_\gep^\infty g(x) \, dx
\\ & =
\int_{\phi(\gep)}^\infty g(x) \, \left(e^{-a x}-1\right) \, dx
+\int_{\phi(\gep)}^1 g(x) \, dx - \int_\gep^1 g(x) \, dx.
\end{split}
\]
Clearly, by assumptions (H1)-(H2) and by dominated convergence
\[
\lim_{\gep\searrow 0} \int_{\phi(\gep)}^\infty g(x) \, \left(e^{-a
x}-1\right) \, dx = \int_0^\infty g(x) \, \left(e^{-a x}-1\right) \,
dx = -\Psi(a).
\]
Now, by assumption (H2)
\[
\begin{split}
& \lim_{\gep\searrow 0} \left[ \int_{\phi(\gep)}^1 g(x) \, dx -
\int_\gep^1 g(x) \, dx \right] = g_0 \, \lim_{\gep\searrow 0}
\left[\int_{\phi(\gep)}^1 \frac1x \, dx - \int_\gep^1 \frac1x \, dx
\right] \\ & = g_0 \, \lim_{\gep\searrow 0} \log
\frac\gep{\phi(\gep)} = -g_0 \, \log \phi'(0).
\end{split}
\]
Then we have obtained \eqref{identity}. \qed

\medskip\noindent
{\it Proof of Proposition \ref{prmain}}. It is enough to apply the
results of Lemma \ref{pr1} and Lemma \ref{est0} to
$\phi(x):=H(s,x)$, $a=\lambda_s$ and $F$ defined in (\ref{defF}).
Positivity of $M_t^{H,\lambda}$ follows from point (3) of Lemma
\ref{pr1}. \qed

\medskip\noindent
{\it Proof of Theorem \ref{main}}. Notice that $M^H=M^{H,\lambda}$
for $\lambda\equiv 0$. By Proposition \ref{prmain}, $M^H$ is a
martingale with expectation 1. Then, for any bounded measurable
$\lambda:\R_+\mapsto\R_+$, by Proposition \ref{prmain} we obtain
\[
\E\left( \exp\left(-\int_0^t \lambda_s \, d \xi^H_s \right) \, M^H_t
\right) = \exp\left(-\int_0^t \Psi(\lambda_s) \, ds \right), \quad
t\geq 0.
\]
The desired result now follows by uniqueness of the Laplace
transform. \qed

\section{Quasi-invariance properties of $(\cL)$-subordinators}
\label{section:special}

In this section we point out two special cases of Theorem
\ref{main}. We consider a measurable function
$k:\R_+\times\R_+\mapsto\R_+$ which satisfies, for some finite
constants $\kappa\geq 1$ and $\alpha\in]0,1[$
\[
|k(s,x)-k(s,y)| \leq \kappa |x-y|^\alpha, \qquad \forall \,
s,x,y\in\R_+,
\]
\[
0< \kappa^{-1} \leq k(s,x) \leq \kappa< \infty, \qquad \forall \,
s,x\in\R_+,
\]
and we set
\[
K(s,x) := \int_0^x k(s,y) \, dy, \qquad \forall \, x,s\geq 0,
\]
Notice that $K(s,\cdot):\R_+\mapsto\R_+$ is necessarily bijective
for all $s\geq 0$.

\subsection{Quasi-invariance of $\xi$ under composition with a diffeomorphism}

Setting
\[
H(s,x) := K(s,\xi_{s-}+x)-K(s,\xi_{s-}), \qquad
h(s,x):=k(s,\xi_{s-}+x), \qquad s,x\geq 0,
\]
we find that
\[
\xi^H_t = \sum_{s\leq t} H(s,\Delta \xi_s) = K(t,\xi_t), \qquad
t\geq 0.
\]
Moreover \eqref{H'1} and \eqref{H'2} are satisfied and Theorem
\ref{main} becomes
\begin{cor}\label{mainc}
The process
\[
G_t^K :=  \exp\left( g_0 \int_0^t \log k(s,\xi_s) \, ds \right)
\prod_{s\leq t} \left[ k(s,\xi_s) \cdot
\frac{g(K(s,\xi_{s})-K(s,\xi_{s-}))}{g(\Delta \xi_s)}\right], \quad
t\geq0,
\]
is a non-negative $(\cF_t)$-martingale with $\E(G_t^K)=1$ and a.s.
$G^K_t>0$. Then we can define a probability measure $\bbP^K$ such
that $\bbP^K_{|\cF_t}=G_t^K\cdot \bbP_{|\cF_t}$ for all $t\geq 0$.
Under $\bbP^K$, $(K(t,\xi_t),t\geq 0)$ is distributed as $(\xi_t,
t\geq 0)$ under $\bbP$.
\end{cor}
\noindent This result can be interpreted by saying that the law of
$(\xi_t)$ is quasi-invariant under (deterministic) non-linear
transformations $(\xi_t, t\geq 0)\mapsto(K(t,\xi_t), t\geq 0)$.

\subsection{Quasi-invariance of $\xi$ under transformations of jumps}

Setting
\[
H(s,x) := K(s,x), \qquad \qquad h(s,x) := k(s,x), \qquad s,x\geq 0,
\]
we find that \eqref{H'1} and \eqref{H'2} are satisfied and Theorem
\ref{main} becomes
\begin{cor}\label{mainc2}
The process
\[
N_t^K :=  \exp\left( g_0 \, \int_0^t \log k(s,0) \, ds\right)
\prod_{s\leq t} \left[k(s,\Delta \xi_s) \cdot \frac{g(K(s,\Delta
\xi_{s}))}{g(\Delta \xi_s)}\right], \quad t\geq 0,
\]
is a non-negative $(\cF_t)$-martingale with $\E(N_t^K)=1$ and a.s.
$N_t^K>0$. Then we can define a probability measure $\bbP^K$ such
that $\bbP^K_{|\cF_t}=N_t^K\cdot \bbP_{|\cF_t}$ for all $t\geq 0$.
Under $\bbP^K$, the process
\[
\xi^K_t = \sum_{s\leq t} K(s,\Delta \xi_s), \qquad t\geq 0,
\]
is distributed as $\xi$ under $\bbP$.
\end{cor}
\noindent This result can be interpreted by saying that the law of
$(\xi_t)$ is quasi-invariant under (deterministic) non-linear
transformation of the jumps of $\xi$: $(\Delta\xi_t, t\geq
0)\mapsto(K(t,\Delta\xi_t), t\geq 0)$.

\subsection{Quasi-invariance properties of the gamma process}

We now write the results of Corollaries \ref{mainc} and \ref{mainc2}
in the special case of the gamma process $(\gamma_t)$. Here
\[
g(x) = \frac{e^{-x}}x, \quad x>0, \qquad g_0=1, \qquad \Psi(\lambda)=\log(1+\lambda).
\]
\begin{cor}\label{prmainga}
We set for all $t\geq 0$
\begin{equation}\label{defMH2}
Y_t^{K} :=  \exp\left( \gamma_t-K(t,\gamma_t) + \int_0^t \log
k(s,\gamma_s)\, ds \right) \prod_{s\leq t} \left[
\frac{k(s,\gamma_s) \cdot \Delta
\gamma_s}{K(s,\gamma_s)-K(s,\gamma_{s-})} \right].
\end{equation}
Then $(Y^K_t)$ is a martingale with $\E(Y^K_t)=1$ and a.s.
$Y_t^K>0$. Hence, we can define a probability measure $\bbP^K$ such
that $\bbP^K_{|\cF_t}=Y_t^K\cdot \bbP_{|\cF_t}$ for all $t\geq 0$.
Under $\bbP^K$, $(K(t,\gamma_t),t\geq 0)$ is distributed as
$(\gamma_t, t\geq 0)$ under $\bbP$.
\end{cor}
\begin{cor}\label{mainc23}
The process
\[
Z_t^K :=  \exp\left( \gamma_t - \sum_{s\leq t} K(s,\Delta\gamma_s)+
\int_0^t \log k(s,0) \, ds\right) \prod_{s\leq t} \left[k(s,\Delta
\gamma_s) \cdot \frac{\Delta \gamma_s}{K(s,\Delta\gamma_s)}\right],
\]
$t\geq 0$, is a non-negative $(\cF_t)$-martingale with $\E(Z_t^K)=1$
and a.s. $Z_t^K>0$. Then we can define a probability measure
$\bbP^K$ such that $\bbP^K_{|\cF_t}=Z_t^K\cdot \bbP_{|\cF_t}$ for
all $t\geq 0$. Under $\bbP^K$, the process
\[
\gamma^K_t = \sum_{s\leq t} K(s,\Delta \gamma_s), \qquad t\geq 0,
\]
is distributed as $(\gamma_t, t\geq 0)$ under $\bbP$.
\end{cor}

\section{Quasi-invariance properties of the Dirichlet Process}
\label{section:dirichlet}

We fix $T>0$ and we denote by $(D^{(T)}_t: t\in[0,T])$ the Dirichlet
process over the time interval $[0,T]$, i.e.
$D^{(T)}_t:=\gamma_t/\gamma_T$ where $(\gamma_t)$ is a gamma
process.  Since $T$ is fixed we omit the superscript $(T)$.

We consider a measurable function $k:[0,T]\times[0,1]\mapsto[0,1]$
which satisfies, for some finite constants $\kappa\geq 1$ and
$\alpha\in]0,1[$
\[
|k(s,x)-k(s,y)| \leq \kappa |x-y|^\alpha, \qquad \forall \,
x,y\in[0,1], \ s\in[0,T],
\]
\[
0< \kappa^{-1} \leq k(s,x) \leq \kappa< \infty, \qquad \forall \,
x\in[0,1], \ s\in[0,T],
\]
and we set
\[
K(s,x) := \int_0^x k(s,y) \, dy, \quad \forall \, x\in[0,1], \
s\in[0,1].
\]

\subsection{Quasi-invariance of $D$ under composition with a diffeomorphism}
We want to give a martingale proof of a relation originally obtained
by von Renesse-Sturm in \cite{vrs}. In this subsection we suppose that $k$ also
satisfies
\[
\int_0^1 k(s,y) \, dy = 1,  \qquad \forall \ s\in[0,T],
\]
so that
\[
K(s,0)=0, \ K(s,1)=1, \qquad \forall \ s\in[0,1].
\]
Notice that $K(s,\cdot):[0,1]\mapsto[0,1]$ is necessarily bijective
for all $s\in[0,T]$. We set for $t<T$
\[
\begin{split}
L_t^{K,T} :=  \left(\frac{1-K(t,D_t)}{1-D_t}\right)^{T-t-1}
\exp\left(\int_0^t \log k(s,D_s) \, ds\right) \prod_{s\leq t}
\left[\frac{k(s,D_s) \cdot \Delta D_s}{K(s,D_s)-K(s,D_{s-})}
\right],
\end{split}
\]
\[
L^{K,T}_T := \frac1{k(T,1)} \, \exp\left(\int_0^T \log k(s,D_s) \,
ds\right) \prod_{s\leq T} \left[\frac{k(s,D_s) \cdot \Delta
D_s}{K(s,D_s)-K(s,D_{s-})}\right].
\]
\begin{theorem}\label{main2}$ $
\begin{enumerate}
\item $(L^{K,T}_t,t\in[0,T])$ is a martingale with respect to the
natural filtration of $D$, such that $\E(L^{K,T}_t)=1$ and a.s.
$L^{K,T}_t>0$, for all $t\in[0,T]$.
\item Under $\bbP^{K,T}:=L^{K,T}_T\cdot \bbP$ the process
$(K(s,D_s), s\in[0,T])$ has the same law as $(D_s, s\in[0,T])$ under
$\bbP$.
%
\end{enumerate}
\end{theorem}
\noindent This theorem gives quasi-invariance of the law of $D$
under non-linear transformations $(D_s, s\in[0,T]) \mapsto
(K(s,D_s), s\in[0,T])$.

\medskip\noindent {\it Proof of Theorem \ref{main2}.} Let first $t<T$.
By the Markov property, for all bounded Borel
$\Phi:\cD([0,t])\mapsto\R_+$
\[
\begin{split}
& \E\left( \Phi(D_s, s\leq t) \right) = \E\left( \Phi(\gamma_s,
s\leq t) \, 1_{(\gamma_t<1)}\, \frac{p_{T-t}(1-\gamma_t)}{p_T(1)}
\right)
\\ & = \E\left( \Phi(\gamma_s, s\leq t)
\, 1_{(\gamma_t<1)}\, (1-\gamma_t)^{T-t-1} \, e^{\gamma_t} \right) \,
\frac{\Gamma(T)}{\Gamma(T-t)}.
\end{split}
\]
Consider the following extension of $K$ to $[0,T]\times\R_+$, that
we still call $K$
\[ K(s,x):= K(s,x)\, 1_{(x\leq 1)} + k(s,1) (x-1)
\, 1_{(x> 1)}, \qquad x\geq 0, \ s\in[0,T].
\]
Let us consider the process $(Y^K_t)$ as defined in \eqref{defMH2}.
Notice that $K(t,\cdot)$ is strictly increasing, so that
$K(t,\gga_t)<1$ iff $\gga_t<1$. Then, for all bounded Borel
$\Phi:{\mathcal D}([0,t])\mapsto\R_+$, $t<T$, by Corollary
\ref{prmainga}
\[
\begin{split}
\E\left( \Phi(K(\cdot,D_\cdot))\, L_t^{K,T} \right) & =
\frac{\Gamma(T)}{\Gamma(T-t)} \, \E\left( 1_{(\gamma_t<1)}\,
\left(1-K(t,\gamma_t)\right)^{T-t-1} \, \Phi(K(\cdot,\gamma_\cdot))
\, e^{K(t,\gamma_t)} \, Y^K_t \right)
\\ & = \frac{\Gamma(T)}{\Gamma(T-t)} \,
\E\left( 1_{(\gamma_t<1)}\, \left(1-\gamma_t\right)^{T-t-1} \,
\Phi(\gamma_\cdot) \, e^{\gamma_t} \right) = \E\left( \Phi(D_\cdot)
\right),
\end{split}
\]
and this concludes the proof for $t<T$.

We consider now the case $t=T$. For all bounded Borel
$\Phi:{\mathcal D}([0,T])\mapsto\R_+$ and $\varphi:\R_+\mapsto\R_+$,
by Corollary \ref{prmainga}
\begin{equation}\label{ccc}
\E\left( \Phi( K(s,\gamma_s), s\in[0,T]) \ \varphi(K(T,\gamma_T)) \
Y^K_T \right) = \E\left( \Phi( \gamma_s, s\in[0,T]) \
\varphi(\gamma_T) \right).
\end{equation}
We set for all $x>0$
\[
Y_T^{K,x} :=  \exp\left( x-K(T,x) + \int_0^T \log k(s,xD_s)\, ds
\right) \prod_{s\leq T} \left[ \frac{k(s,xD_s) \cdot x\Delta
D_s}{K(s,xD_s)-K(s,xD_{s-})} \right].
\]
In the right hand side of \eqref{ccc} we condition on the value of
$\gamma_T$, obtaining
\[
\E\left( \Phi( \gamma_\cdot) \, \varphi(\gamma_T) \right) =
\int_0^\infty p_T(y) \, \E\left( \Phi( yD_\cdot)  \right) \varphi(y)
\, dy.
\]
In the left hand side of \eqref{ccc}, conditioning on the value of
$\gamma_T$, we obtain
\[
\E\left( \Phi( K(\cdot,\gamma_\cdot)) \ \varphi(K(T,\gamma_T)) \,
Y^K_T \right) = \int_0^\infty p_T(x) \, \E\left( \Phi(
K(\cdot,xD_\cdot)) \, Y^{K,x}_T \right) \varphi(K(T,x)) \, dx.
\]
In order to compare this expression with the one above for the right
hand side, we use the change of variable $x=K(T,y)$. To this aim, we
denote by $C:\R_+\mapsto\R_+$ the inverse of $K(T,\cdot)$, i.e. we
suppose that $K(T,C(x))=x$ for all $x\geq 0$. Then we have
\[
\begin{split}
& \E\left( \Phi( K(\cdot,\gamma_\cdot)) \ \varphi(K(T,\gamma_T)) \,
Y^K_T \right)  = \int_0^\infty p_T(x) \, \E\left( \Phi(
K(\cdot,xD_\cdot)) \, Y^{K,x}_T \right) \varphi(K(T,x)) \, dx
 \\ & = \int_0^\infty p_T(C(y)) \, \E\left( \Phi(
K(\cdot,C(y)D_\cdot)) \, Y^{K,C(y)}_T \right) \varphi(y) \, C'(y) \,
dy.
\end{split}
\]
Since this is true for any bounded measurable
$\varphi:\R_+\mapsto\R_+$, we obtain for all $y>0$
\[
\frac{p_T(C(y)) \, C'(y)}{p_T(y)} \, \E\left( \Phi(
K(\cdot,C(y)D_\cdot)) \ Y^{K,C(y)}_T \right) = \E\left( \Phi(
yD_\cdot) \right).
\]
For $y=1$, since $K(T,1)=1=C(1)$, we obtain the desired result
\[
\E\left( \Phi( K(\cdot,D_\cdot)) \ L^{K,T}_T \right) = \E\left(
\Phi( D_\cdot) \right). \qed
\]

\begin{rem}{\rm Von Renesse-Sturm prove the second result of Theorem
\ref{main2} in  \cite{vrs}. The proof there hinges on explicit
computations related to the finite-dimensional distributions of $D$.
}
\end{rem}

\subsection{Quasi-invariance of $D$ under transformation of the jumps}

Again, we consider the Dirichlet process $(D^{(T)}_t, t\in[0,T])$,
and we drop the superscript $(T)$, since $T$ is fixed. We set
\[
\Delta D_s := D_s-D_{s-}, \qquad D^K_t:= \frac{\sum_{s\leq t}
K(s,\Delta D_s)}{\sum_{s\leq T} K(s,\Delta D_s)}, \qquad t\in[0,T].
\]
\begin{theorem}\label{main3}
The laws of $(D^K_t, t\in[0,T])$ and $(D_t, t\in[0,T])$ are
equivalent.
\end{theorem}
\noindent In the proof of Theorem \ref{main3} we also compute
explicitly the Radon-Nikodym density. Handa \cite{handa} considers
the particular case $K(s,x)=c(s) \, x$, where $c:[0,T]\mapsto\R_+$
is measurable.

\medskip\noindent {\it Proof}. We set
\[
\gga^K_t:= \sum_{s\leq t} K(s,\Delta\gga_s), \qquad t\geq 0.
\]
Since $(D_t,t\in[0,T])$ is a gamma bridge, then the law of $(D^K_t,
t\in[0,T])$ coincides with the law of $(\gamma^K_t/\gamma^K_T,
t\in[0,T])$ under the conditioning $\{\gamma_T=1\}$.

We define $J:[0,T]\times\R_+\mapsto\R$, such that, for all
$s\in[0,T]$, $J(s,K(s,x))=x$ for all $x\geq 0$. In other words,
$J(s,\cdot)$ is the inverse of $K(s,\cdot)$. In particular, notice
that
\[
\left(\gamma^K\right)^J = \gamma.
\]
By Corollary \ref{mainc23}, for all $\Phi:\cD([0,T])\mapsto\R$
bounded and Borel
\begin{equation}\label{coco}
\E\left(\Phi\left( \frac{\gamma^K_s}{\gamma^K_T}, s\leq T \right) \,
\varphi(\gga_T) \, Z^K_T \right) = \E\left(\Phi\left(
\frac{\gamma_s}{\gamma_T}, s\leq T \right) \, \varphi(\gga_T^J)
\right).
\end{equation}
Notice that
\[
\gga_T^J = \sum_{s\leq T} J(s,\Delta\gga_s) = \sum_{s\leq T}
J\left(s,\gga_T\cdot\Delta D_s\right) = \psi_D( \gga_T),
\]
where $D_t:=\gamma_t/{\gamma_T}$, $t\in[0,T]$, is independent of
$\gamma_T$ and
\[
\psi_D(x) := \sum_{s\leq T} J\left(s,x\cdot\Delta D_s\right), \qquad
x\geq 0.
\]
Notice that $\psi_D:\R_+\mapsto\R_+$ is $C^1$ and by dominated
convergence
\[
\psi_D'(x) =  \sum_{s\leq T} \Delta D_s\cdot J\left(s,x\cdot\Delta
D_s\right) \geq \kappa^{-1}>0, \qquad \forall \ x\geq 0,
\]
since $\Delta D_s\geq 0$ and $\sum_{s\leq T} \Delta D_s = 1$. Also
by dominated convergence, $\psi_D'$ is continuous. Therefore
$\psi_D:\R_+\mapsto\R_+$ is invertible, with $C^1$ inverse
$\zeta_D:=\psi_D^{-1}$. In the sequel, We may write $\zeta_D^K$ for
$\zeta_D$, in order to stress that it also depends on $K$. Then, for
all $\varphi:\R_+\mapsto\R$ bounded and Borel, we obtain by
\eqref{coco}
\begin{eqnarray}\label{eqna1}
\nonumber & & \E\left( \Phi\left( \frac{\gamma^K_s}{\gamma^K_T},
s\leq T \right) \, \varphi(\gga_T) \, Z^K_T \right) =
\E\left(\Phi\left( D_s, s\leq T \right) \, \varphi(\psi_D( \gga_T))
\right) \\ \nonumber & & =  \E\left(\Phi\left( D_s, s\leq T \right)
\int_0^\infty p_T(y) \, \varphi(\psi_D(y)) \, dy \right) = \big[ x=
\psi_D(y) \big]
\\ & & = \int_0^\infty \varphi(x) \, \E\Big(\Phi\left( D_s, s\leq T \right) \, p_T(\zeta_D(x))
\, \zeta_D'(x) \, \Big) \, dx.
\end{eqnarray}
Now, setting for all $t\in[0,T]$
\[
D^{K,x}_t:= \frac{\sum_{s\leq t} K(s,x\cdot\Delta D_s)}{\sum_{s\leq
T} K(s,x\cdot\Delta D_s)},
\]
\[
U^{K,x}_T := \exp\left( x - \sum_{s\leq T} K(s,x\Delta D_s)+
\int_0^T \log k(s,0) \, ds\right) \prod_{s\leq T}
\left[\frac{k(s,x\Delta D_s) \cdot x\Delta D_s}{K(s,x\Delta
D_s)}\right],
\]
then we have
\begin{equation}\label{eqna2}
\E\left( \Phi\left( \frac{\gamma^K_s}{\gamma^K_T}, s\leq T \right)
\, \varphi(\gga_T) \, Z^K_T \right) = \int_0^\infty \varphi(x) \,
\E\left( \Phi\left( D^{K,x}_s, s\leq T \right) \cdot U^{K,x}_T
\right) \, p_T(x) \, dx.
\end{equation}
Since $D^{K,1}=D^{K}$, setting
\[
U^K_T := U^{K,1}_T = \exp\left(1- \sum_{s\leq T} K(s,\Delta D_s)
+\int_0^T \log k(s,0)\, ds \right) \prod_{s\leq T} \left[
\frac{k(s,\Delta D_s) \cdot \Delta D_{s}}{K(s,\Delta D_s)}\right],
\]
we obtain by \eqref{eqna1} and \eqref{eqna2} for $x=1$
\[
\E\left( \Phi\left( D^K_s, s\leq T \right) \cdot U^K_T \right) \, =
\E\left(\Phi\left( D_s, s\leq T \right) \,
\frac{p_T(\zeta_D^K(1))}{p_T(1)} \, (\zeta_D^K)'(1) \right). \qed
\]

\section{Stochastic differential equations driven by $(\cL)$-subordinators}
\label{sdes}

In this section we give an application to stochastic differential
equations driven by a $(\cL)$-subordinator $\xi$. See \cite{bass}
for a survey of SDEs driven by L\'evy processes.

We consider the SDE
\begin{equation}\label{sde}
dX_t = m(t,X_{t-}) \, d\xi_t, \qquad X_0=0,
\end{equation}
where
\begin{enumerate}
\item $m:\R_+\times\R_+\mapsto(0,+\infty)$ is measurable;
\item $m$ and $1/m$ are bounded
\item $\R_+\ni a\mapsto m(s,a)$ is Lipschitz, uniformly in $s\geq
0$.
\end{enumerate}
Then we have
\begin{theorem}\label{thsde}
There exists a pathwise-unique solution of \eqref{sde} and the law of $(X,\xi)$
under $\bbP$ coincides with the law of $(\xi,\xi^H)$ under $\bbP^H$, where
\begin{equation}\label{h}
H(s,x) := \frac x{m(s,\xi_{s-})}, \qquad s\geq 0, \ x\geq 0.
\end{equation}
\end{theorem}
\noindent {\it Proof}. Let $T>0$ and denote by $\cI([0,T])$ the set
of all bounded increasing functions $\omega:[0,T]\mapsto \R_+$. We
define the map $\Lambda_T:\cI([0,T])\mapsto \cI([0,T])$
\[
\Lambda_T(\omega)(t):=\int_0^t m(s,\omega_{s-}) \, d\xi_s, \qquad
t\in[0,T].
\]
For $L$ large enough, $\Lambda_T$ is a contraction in $\cI([0,T])$ with respect to the metric
\[
d_L(\omega,\omega') := \sup_{t\in[0,T]} \, e^{-Lt} \, |\omega_t-\omega_t'|,
\]
and the solution of \eqref{sde} on the time
interval $[0,T]$ is the unique fixed point $X$ of $\Lambda_T$.
Moreover, there exists a measurable map $W_T:\cI([0,T])\mapsto
\cI([0,T])$, such that $X=W_T(\xi_{|[0,T]})$.

Let us define $H$ as in \eqref{h}, and set $\xi^H$ as in
\eqref{defxih}
\[
\xi^H_t := \sum_{s\leq t} H(s,\Delta \xi_s) = \int_0^t \frac
1{m(s,\xi_{s-})} \, d\xi_s
\]
Note that
\[
d\xi^H_t = \frac 1{m(t,\xi_{t-})} \, d\xi_t \ \Longrightarrow \
d\xi_t = {m(t,\xi_{t-})} \, d\xi^H_t.
\]
Then, $\xi_{|[0,T]}=W_T(\xi^H_{|[0,T]})$ for any $T>0$. On the
other hand, by Theorem \ref{main}, $\xi^H$ under $\bbP^H$ has the same
law as $\xi$ under $\bbP$, and this concludes the proof.
\qed


\begin{thebibliography}{abc99xys}

\bibitem{apple} D. Applebaum (2004), {\it L\'evy processes and stochastic calculus},
Cambridge Studies in Advanced Mathematics, {\bf 93}, Cambridge
University Press.

\bibitem{bass} R.F. Bass (2004), {\it Stochastic differential equations with jumps},
Probab. Surv.  1, 1--19.

\bibitem{baxter} M. Baxter (2006), {\it L\'evy process dynamic modelling of singlename
credits and CDO tranches}, Working Paper, Nomura Fixed Income Quant
Group.

\bibitem{bertoin} J. Bertoin (1997), {\it Subordinators: examples and applications},
in: {Lectures on probability theory and statistics
(Saint-Flour,1997)}, Springer Verlag, Lecture Notes in Math. 1717.

\bibitem{handa} K. Handa (2002), {\it Quasi-invariance and reversibility in the Fleming-Viot process},
Probab. Theory Related Fields  {\bf 122},  no. 4, 545--566.

\bibitem{vrs} M.-K. von Renesse, K.-T. Sturm (2007), {\it Entropic Measure and Wasserstein Diffusion},
preprint, posted on http://arxiv.org/abs/0704.0704.

\bibitem{reyo} D. Revuz, and M. Yor (1991),
{\it Continuous Martingales and Brownian Motion}, Springer Verlag.

\bibitem{sato} K.I. Sato (1999), {\it L\'evy processes and infinitely divisible distributions},
Cambridge Studies in Advanced Mathematics, {\bf 68}, Cambridge
University Press, Cambridge.

\bibitem{tvy} N. Tsilevich, A. Vershik, M. Yor (2001), {\it An infinite-dimensional analogue of the Lebesgue
measure and distinguished properties of the gamma process},  J.
Funct. Anal.  {\bf 185},  no. 1, 274--296.

\bibitem{tvy2} N. Tsilevich, A. Vershik, M. Yor (2004), {\it On the Markov-Krein identity and the
quasi-invariance of the gamma process},    J. Math. Sci. {\bf 121},
no. 3, 2303--2310.

\bibitem{yor} M. Yor (2007), {\it Some remarkable properties of Gamma processes},
in: Advances in Mathematical Finance, Festschrift volume in honour
of Dilip Madan, eds. R.J. Elliott, M.C. Fu, R.A. Jarrow, and J.-Y.J.
Yen, Birkh\"auser/Springer.




\end{thebibliography}
\end{document}